\documentclass[12pt]{article}
\usepackage{amsthm}
\usepackage{amsmath,amsfonts,amssymb}
\usepackage{color,graphics}
\usepackage{indentfirst}
\usepackage{graphicx}
\usepackage{geometry}
\usepackage{titlesec}
\usepackage{caption}
\usepackage{authblk}
\usepackage[misc]{ifsym}  
\usepackage{skak} 
\newtheorem{dingyi}{Definition}[section]
\newtheorem{yinli}[dingyi]{Lemma}
\newtheorem{dingli}[dingyi]{Theorem}

\newtheorem{mingti}[dingyi]{Proposition}

\newcommand{\GL}{\mathrm{GL}}

\newcommand{\spa}{\mathrm{span}}
\newcommand{\bol}{Bollob\'as}
\newcommand{\lov}{Lov\'asz}
\newcommand{\fu}{F\"{u}redi}
\newcommand{\ekr}{Erd\"{o}s-Ko-Rado}
\newcommand{\heg}{Heged\"{u}s}

\newcommand{\ELTE}{E\"{o}tv\"{o}s Lor\'{a}nd University}
\newcommand{\ELTEAdress}{P\'{a}zm\'{a}ny P\'{e}ter s\'{e}t\'{a}ny 1/C, Budapest, Hungary, H-1117}

\title{A~\bol-type theorem on singular linear spaces}

\author[$\ast\spadesuit$]{Erfei Yue}
\author[$\heartsuit$]{Benjian Lv}
\author[$\diamondsuit$]{P\'{e}ter Sziklai}
\author[$\clubsuit$]{Kaishun Wang}

\affil[$\spadesuit$]{\footnotesize Institute of Mathematics, \ELTE,

\ELTEAdress,

\Letter\ yef9262@mail.bnu.edu.cn}

\affil[$\heartsuit$]{\footnotesize Laboratory of Mathematics and Complex Systems (Ministry of Education),

School of Mathematical Sciences, Beijing Normal University, Beijing 100875, China

\Letter\ bjlv@bnu.edu.cn}
\affil[$\diamondsuit$]{\footnotesize Institute of Mathematics, \ELTE,

\ELTEAdress,

\Letter\ peter.sziklai@ttk.elte.hu}
\affil[$\clubsuit$]{\footnotesize Laboratory of Mathematics and Complex Systems (Ministry of Education),

School of Mathematical Sciences, Beijing Normal University, Beijing 100875, China

\Letter\ wangks@bnu.edu.cn}

\date{}

\begin{document}
\maketitle

\begin{center}
\textbf{Abstract}
\end{center}
\bol-type theorem determines the maximum cardinality of a~\bol~system of sets.
The original result has been extended to various mathematical structures beyond sets, including vector spaces and affine spaces.
This paper generalizes the \bol-type theorem to singular linear spaces, and determine the maximum cardinality of (skew)~\bol~systems on them.

\section{Introduction}

Let~$[n]=\{1,\ldots,n\}$ be a set of cardinality~$n$.
The following concepts of~\bol~systems and skew~\bol~systems (of sets) were introduced by~\bol~\cite{Sets} and Frankl~\cite{Skew} respectively.

\begin{dingyi}
Suppose~$\mathcal{P}=\{(A_i,B_i)\mid i\in[m]\}$ is a family of pairs of sets, where~$A_i,B_i\subseteq [n]$, and~$A_i\cap B_i=\emptyset$.
Then~$\mathcal{P}$ is called a~\emph{\bol~system} if~$A_i\cap B_j\neq\emptyset,\forall i\neq j$,
and a \emph{skew~\bol~system} if~$A_i\cap B_j\neq\emptyset,\forall i<j$.
\end{dingyi}

The following theorem showing the maximum cardinality of a~\bol~system was first proved by~\bol~\cite{Sets} in 1965,
and later proved by Jaeger and Payan~\cite{FRA}, Katona~\cite{Katona}, and Tarj\'an~\cite{Tarjan} independently.
In 1982, Frankl~\cite{Skew} prove that the result remains true even when its condition of~\bol~system are replaced with skew~\bol~system.
This generalization is useful in automata theory~\cite{Pin}.

\begin{dingli}[\bol~\cite{Sets},~Frankl~\cite{Skew}]\label{Th:uniform}
Let~$\mathcal{P}=\{(A_i,B_i)\mid i\in[m]\}$ be a (skew)~\bol~system, where~$A_i,B_i\subseteq [n]$, and~$|A_i|=a,|B_i|=b$ for any~$i$.
Then~$m\leqslant\binom{a+b}{a}$.
\end{dingli}

In 1985, Alon~\cite{Alon} proved the following theorem, which is a variation of~\bol-type theorem,
and determine the maximum cardinality of a (skew)~\bol~system on~$r$-partitions of sets.

\begin{dingli}[Alon~\cite{Alon}]\label{Th:Sets}
Suppose~$X$ is a disjoint union of some sets~$X_1,\ldots,X_r$,~$|X_k|=n_k$, and~$n_1+\cdots+n_r=n$.
A skew~\bol~system of subsets~$\mathcal{P}=\{(A_i,B_i)\mid i\in[m]\}$ of~$X$ satisfies that
\begin{equation*}
|A_i|=a,|B_i|=b,|A_i\cap X_k|=a_k,|B_i\cap X_k|=b_k, \forall i\in[m],k\in[r].
\end{equation*}
Then we have
\begin{equation*}
m\leqslant\prod_{k=1}^r\binom{a_k+b_k}{a_k}.
\end{equation*}
\end{dingli}

As an analog, the~\bol~systems and skew~\bol~systems of vector spaces are defined as follows.

\begin{dingyi}
Let~$\mathcal{P}=\{(A_i,B_i)\mid i\in[m]\}$ be a family of pairs of subspaces of a fixed vector space~$V$, where~$\dim(A_i\cap B_i)=0$.
We say that~$\mathcal{P}$ is a\emph{~\bol~system} if~$\dim(A_i\cap B_j)>0,\forall i\neq j$,
and a \emph{skew~\bol~system} if~$\dim(A_i\cap B_j)>0,\forall i<j$.
\end{dingyi}

In 1977, leveraging powerful algebraic techniques,~\lov~\cite{LinearSpaces,Lov1,Lov2} extended Theorem~\ref{Th:uniform} to the realm of matroids.
This remarkable extension opened new avenues for solving combinatorial problems using the exterior product method,
and solidified the~\bol-type theorem as one of the central problems in extremal set theory.
Maybe for the sake of simplicity,~\fu~\cite{Threshold} reformulated it as follows for real vector spaces.

\begin{dingli}[\lov~\cite{LinearSpaces,Lov1,Lov2}]\label{Th:Lovasz}
Let~$\mathcal{P}=\{(A_i,B_i)\mid i\in[m]\}$ be a (skew)~\bol~system of subspaces of~$\mathbb{R}^n$.
If~$\dim(A_i)=a,\dim(B_i)=b$ for all~$i$, then~$m\leqslant\binom{a+b}{a}$.
\end{dingli}

The conclusion of Theorem~\ref{Th:Lovasz} is true for vector spaces over an arbitrary field.
For the case of~$\mathbb{Q}^n$ or~$\mathbb{C}^n$, even the proofs keep the same. For other fields such as finite fields,
some technologies need to be adopted. More details about this can be found in~\cite{Book}.
After~\fu, other scholars~(such as Scott and Wilmer~\cite{Scott}, and Yu, Kong, Xi, Zhang, and Ge~\cite{HemiBundled}) tend to express their theorems in the same way,
even though they hold more generally.

\bol-type theorem and its variations have attracted significant interest among mathematicians.
In 1984,~\fu~\cite{Threshold} established threshold (or~$t$-intersecting) version of~\bol-type theorem (for both sets and spaces).
This work paved the way for further extensions by Zhu~\cite{tl}, Talbot~\cite{Inequality}, and Kang, Kim, and Kim~\cite{Inequality2}.
Recent research has explored even more variations of~\bol-type theorem,
such as~\cite{Alon} on partitions of sets,~\cite{AffineSpaces} on affine spaces,~\cite{weakly} on weakly~\bol~systems,~\cite{kTuples} on~$k$-tuples,
~\cite{Scott} on nonuniform~\bol~systems,~\cite{Hegedus} on nonuniform skew~\bol~systems, and~\cite{HemiBundled} on hemi-bundled two families.

Fix an~$(n+l)$-dimensional vector space~$V$, and an~$l$-dimensional subspace~$W$.
The stabilizer~$G$ of~$W$ for the general linear group~$\GL(V)$ is called the \emph{singular linear group}.
The space~$V$ together with the action of~$G$ is called a \emph{singular linear space}.
Consider the action of~$G$ on the subspaces of~$V$. It is easy to see that subspaces~$P$ and~$Q$ fall in the same orbit if and only if
\begin{equation*}
\dim(P)=\dim(Q),\dim(P\cap W)=\dim(Q\cap W).
\end{equation*}
Hence we say that an~$s$-dimensional subspace~$P$ is of type~$(s,t)$, if~$\dim(P\cap W)=t$.
The concept of singular linear spaces over a finite field was introduced in~\cite{AssociationSchemes,SingularLinearSpace},
and has a wide range of applications in many fields such as association scheme and coding theory.

In this paper, we generalize the~\bol-type theorem to singular linear spaces.
More specifically, for fixed parameters~$a,a_0,b,b_0,n,l$ and an~$l$-dimensional subspace~$W$ of~$V=\mathbb{R}^{n+l}$,
suppose~$\mathcal{P}=\{(A_i,B_i)\mid i\in[m]\}$ is a (skew)~\bol~system, where~$A_i$'s and~$B_i$'s are~$(a,a_0)$-type and~$(b,b_0)$-type subspaces of~$V$ correspondingly.
We will explore~$M(a,a_0,b,b_0,n,l)$, the maximum cardinality of~$\mathcal{P}$.
The following theorem is our main result, which is a generalization of Theorem~\ref{Th:Lovasz}.
Actually, Theorem~\ref{Th:Lovasz} is a special case of Theorem~\ref{Th:Spaces} for~$l=a_0=b_0=0$.
Meanwhile, Theorem~\ref{Th:Spaces} is an extension of Theorem~\ref{Th:Sets} for the case~$r=2$.

\begin{dingli}\label{Th:Spaces}
Let~$V=\mathbb{R}^{n+l}$, and~$W$ be a fixed~$l$-dimensional subspace of~$V$.
Suppose that~$\mathcal{P}=\{(A_i,B_i)\mid i\in[m]\}$ is a skew~\bol~system of subspaces,
where~$A_i$'s are~$(a,a_0)$-type subspaces and~$B_i$'s are~$(b,b_0)$-type subspaces. If~$l\leqslant\min\{a_0+b,a+b_0\}$, then we have
\begin{equation}\label{Eq:main}
m\leqslant\sum_{k=0}^{l-a_0-b_0}\binom{l}{a_0+k}\binom{a+b-l}{a-a_0-k}.
\end{equation}
\end{dingli}

This paper is organized as follows. In Section 2, we make a brief introduction of exterior products and general position arguments,
which will play important roles in the proof. Section 3 presents a detailed proof of Theorem~\ref{Th:Spaces}.
In Section 4, we construct two examples. One of them shows the tightness of Inequality~(\ref{Eq:main}) in Theorem~\ref{Th:Spaces},
and the other shows the condition~$l\leqslant\min\{a_0+b,a+b_0\}$ is not redundant,
by explores the lower bound of~$M(a,a_0,b,b_0,n,l)$ when~$l$ is large.

\section{Preliminaries}

A handy tool to deal with (skew)~\bol~systems of spaces is exterior algebra.
Suppose that~$\underline{\alpha_1},\ldots,\underline{\alpha_k}$ are~$n$-dimensional (raw) vectors in~$\mathbb{R}^n$, let
\begin{equation*}
A=\begin{pmatrix}
\underline{\alpha_1} \\
\vdots \\
\underline{\alpha_k}
\end{pmatrix}
\end{equation*}
be a~$k\times n$ matrix. For~$I\in\binom{[n]}{k}$, denote~$A_I$ the~$k\times k$ submatrix of~$A$ consisting of the~$k$ columns labeled by~$I$.
Suppose~$\binom{[n]}{k}=\{I_1,\ldots,I_m\}$, where~$m=\binom{n}{k}$.
Then the exterior product (or wedge product) of~$\underline{\alpha_1},\ldots,\underline{\alpha_k}$ is an~$\binom{n}{k}$-dimensional vector
\begin{equation*}
\underline{\alpha_1}\wedge\cdots\wedge\underline{\alpha_k}:=(|A_{I_1}|,\ldots,|A_{I_m}|).
\end{equation*}
Note that if~$n=k$ then the exterior product is just the determinant of~$A$.
Suppose~$\underline{\alpha}=\underline{\alpha_1}\wedge\cdots\wedge\underline{\alpha_s}$
and~$\underline{\beta}=\underline{\beta_1}\wedge\cdots\wedge\underline{\beta_t}$, denote
\begin{equation*}
\underline{\alpha}\wedge\underline{\beta}=\underline{\alpha_1}\wedge\cdots\wedge\underline{\alpha_s}
\wedge\underline{\beta_1}\wedge\cdots\wedge\underline{\beta_t}
\end{equation*}
for the sake of simplicity.
The above simplified version of the definition of exterior product comes from~\cite{Book}, in where more details and propositions can be found.
For the original abstract definition of exterior product and exterior algebra, one may refer to any textbooks of multilinear algebra.
Here we emphasize the following two lemmas, which are basic conclusions in multilinear algebra, but play a vital role in the proof of Theorem~\ref{Th:Spaces}.

\begin{yinli}
Let~$V=\mathbb{R}^n$, and~$V^k$ be the Cartesian product of~$k$ many of~$V$. Then
\begin{equation*}
f:V^k\rightarrow \mathbb{R}^{\binom{n}{k}},
\quad (\underline{\alpha_1},\ldots,\underline{\alpha_k})\mapsto\underline{\alpha_1}\wedge\cdots\wedge\underline{\alpha_k}
\end{equation*}
is a~$k$-linear mapping.
\end{yinli}

\begin{yinli}
Suppose that~$\underline{\alpha_1},\ldots,\underline{\alpha_k}$ are~$k$ many of~$n$-dimensional vectors.
Then they are linearly independent if and only if~$\underline{\alpha_1}\wedge\cdots\wedge\underline{\alpha_k}\neq\underline{0}$.
\end{yinli}

To prove Theorem~\ref{Th:Spaces}, we also need the following concept and lemma for general position argument.
Suppose~$V$ and~$W$ are~$n$-dimensional and~$k$-dimensional vector spaces over field~$F$,~$U_1,\cdots,U_m$ are subspaces of~$V$, and~$\dim(U_i)=r_i$.
We say linear mapping~$\phi:V\rightarrow W$  is in general position with~$U_i$'s, if
\begin{equation*}
\dim(\phi(U_i))=\min\{r_i,k\},\quad i=1,\ldots,m.
\end{equation*}
Note that for every subspace~$U_i$,~$\min\{r_i,k\}$ is the maximum possible dimension of~$\phi(U_i)$.
The following lemma shows that for fixed~$V$,~$W$, and~$U_i$'s, a linear mapping~$\phi:V\rightarrow W$ that in general position with~$U_i$'s always exists,
as long as~$|F|$ is large enough.

\begin{yinli}[\cite{Book}]\label{GP3}
Suppose~$V$ and~$W$ are~$n$-dimensional and~$k$-dimensional vector spaces over field~$F$,~$U_1,\ldots,U_m$ are subspaces of~$V$, and~$\dim(U_i)=r_i$.
If~$|F|>(n-k)(m+1)$ then there exists a linear mapping~$\phi:V\rightarrow W$ that is in general position with~$U_i$'s.
\end{yinli}

\section{Proof of the main theorem}

We first consider about a spacial case, and then the general case.

\textbf{Case 1:}~$n+l=a+b$.

Pick~$\underline{\epsilon_1},\ldots,\underline{\epsilon_l}$, a basis of~$W$,
and expand it to~$\underline{\epsilon_1},\ldots,\underline{\epsilon_l},\underline{\gamma_1},\ldots,\underline{\gamma_n}$, a basis of~$V$.
For any~$i$, let~$\underline{\alpha^{(i)}_1},\ldots,\underline{\alpha^{(i)}_{a_0}}$ and
~$\underline{\beta^{(i)}_1},\ldots,\underline{\beta^{(i)}_{b_0}}$ be bases of~$A_i\cap W$ and~$B_i\cap W$.
Expand them to~$\underline{\alpha^{(i)}_1},\ldots,\underline{\alpha^{(i)}_a}$ and~$\underline{\beta^{(i)}_1},\ldots,\underline{\beta^{(i)}_b}$,
bases of~$A_i$ and~$B_i$ accordingly.
Note that~$A_i\cap B_i=\{0\}$, so~$\underline{\alpha^{(i)}_1},\ldots,\underline{\alpha^{(i)}_{a_0}},
\underline{\beta^{(i)}_1},\ldots,\underline{\beta^{(i)}_{b_0}}$ is a linearly independent system in~$W$ for any~$i$.
Again, expand it to
\begin{equation*}
\left\{\underline{\alpha^{(i)}_1},\ldots,\underline{\alpha^{(i)}_{a_0}},\underline{\beta^{(i)}_1},\ldots,\underline{\beta^{(i)}_{b_0}},\underline{\omega^{(i)}_1},
\ldots,\underline{\omega^{(i)}_{l-a_0-b_0}}\right\},
\end{equation*}
a basis of~$W$.
Then\begin{equation*}
S_i:=\left\{\underline{\alpha^{(i)}_1},\ldots,\underline{\alpha^{(i)}_{a_0}},\underline{\beta^{(i)}_1},\ldots,\underline{\beta^{(i)}_{b_0}},\underline{\omega^{(i)}_1},
\ldots,\underline{\omega^{(i)}_{l-a_0-b_0}},\underline{\gamma_1},\ldots,\underline{\gamma_n}\right\}
\end{equation*}
forms a basis of~$V$. Let
\begin{equation*}
M_i=\{\underline{\beta^{(i)}_1},\ldots,\underline{\beta^{(i)}_{b_0}},\underline{\omega^{(i)}_1},\ldots,\underline{\omega^{(i)}_{l-a_0-b_0}}\},
\quad N_i=M_i\cup\{\underline{\gamma_1},\ldots,\underline{\gamma_n}\}.
\end{equation*}

For any~$i$, let
\begin{equation*}
\begin{array}{l}
\underline{u_i}=\underline{\alpha^{(i)}_1}\wedge\cdots\wedge\underline{\alpha^{(i)}_{a_0}}, \\
\underline{u'_i}=\underline{\alpha^{(i)}_{a_0+1}}\wedge\cdots\wedge\underline{\alpha^{(i)}_a}, \\
\underline{\bar{u}_i}=\underline{u_i}\wedge \underline{u'_i},
\end{array}
\end{equation*}
and
\begin{equation*}
\begin{array}{l}
\underline{v_i}=\underline{\beta^{(i)}_1}\wedge\cdots\wedge\underline{\beta^{(i)}_{b_0}}, \\
\underline{v'_i}=\underline{\beta^{(i)}_{b_0+1}}\wedge\cdots\wedge\underline{\beta^{(i)}_b}, \\
\underline{\bar{v}_i}=\underline{v_i}\wedge \underline{v'_i}.
\end{array}
\end{equation*}
Due to the definition of skew~\bol~system, we have~$\underline{\bar{u}_i}\wedge\underline{\bar{v}_i}\neq 0$,
and~$\underline{\bar{u}_i}\wedge\underline{\bar{v}_j}=0$ if~$i<j$.

For a fixed~$i$, since~$S_i$ is a basis of~$V$, we can express~$\underline{u'_i}$ in the form of
\begin{equation*}
\sum_h D_h\underline{\xi_{h,1}}\wedge\cdots\wedge\underline{\xi_{h,a-a_0}},
\end{equation*}
where~$\underline{\xi_{h,k}}\in S_i$, and~$D_h$ is a real number.
Now we classify the terms in the following way. We say~$D_h\underline{\xi_{h,1}}\wedge\cdots\wedge\underline{\xi_{h,a-a_0}}$ is in class
\begin{equation*}
\begin{array}{ll}
C_1: & \textrm{if~} \{\underline{\xi_{h,1}},\ldots,\underline{\xi_{h,a-a_0}}\}\subseteq N_i,
~\textrm{and~}|\{\underline{\xi_{h,1}},\ldots,\underline{\xi_{h,a-a_0}}\}\cap M_i|\leqslant l-a_0-b_0, \\
C_2: & \textrm{if~} \{\underline{\xi_{h,1}},\ldots,\underline{\xi_{h,a-a_0}}\}\subseteq N_i,
~\textrm{and~}|\{\underline{\xi_{h,1}},\ldots,\underline{\xi_{h,a-a_0}}\}\cap M_i|> l-a_0-b_0, \\
C_3: & \textrm{if~} \{\underline{\xi_{h,1}},\ldots,\underline{\xi_{h,a-a_0}}\}\not\subseteq N_i.
\end{array}
\end{equation*}
Let
\begin{equation*}
\underline{u'_i}=\underline{u'_i(1)}+\underline{u'_i(2)}+\underline{u'_i(3)},
\end{equation*}
where
\begin{equation*}
\underline{u'_i(d)}=\sum\limits_{C_d}D_h\underline{\xi_{h,1}}\wedge\cdots\wedge\underline{\xi_{h,a-a_0}},\quad d=1,2,3.
\end{equation*}
Then
\begin{equation*}
\underline{\bar{u}_i}=\underline{u_i}\wedge \underline{u'_i}=\underline{u_i}\wedge(\underline{u'_i(1)}+\underline{u'_i(2)}+\underline{u'_i(3)})
=\underline{u_i}\wedge \underline{u'_i(1)}+\underline{u_i}\wedge \underline{u'_i(2)}.
\end{equation*}
The last term vanished due to the proposition of exterior product. Let~$\underline{\tilde{u}_i}=\underline{u_i}\wedge \underline{u'_i(1)}$.
Note that
\begin{equation*}
0\neq\underline{\bar{u}_i}\wedge\underline{\bar{v}_i}=
\underline{u_i}\wedge \underline{u'_i(1)}\wedge\underline{\bar{v}_i}+\underline{u_i}\wedge \underline{u'_i(2)}\wedge\underline{\bar{v}_i},
\end{equation*}
but~$\underline{u'_i(2)}\wedge \underline{v_i}=\underline{0}$, so~$\underline{\tilde{u}_i}\wedge\underline{\bar{v}_i}\neq 0$.
On the other hand, if~$i<j$ then
\begin{equation*}
0=\underline{\bar{u}_i}\wedge\underline{\bar{v}_j}=
\underline{u_i}\wedge \underline{u'_i(1)}\wedge\underline{\bar{v}_j}+\underline{u_i}\wedge \underline{u'_i(2)}\wedge\underline{\bar{v}_j}.
\end{equation*}
Note that every term in~$\underline{u'_i(2)}$ contains more than~$l-a_0-b_0$ members in~$W$.
Then the vectors~$\underline{\alpha^{(i)}_1},\ldots,\underline{\alpha^{(i)}_{a_0}},\underline{\beta^{(j)}_1},\ldots,\underline{\beta^{(j)}_{b_0}}$
and entries of~$\underline{u'_i(2)}$ are linearly dependent, so~$\underline{u_i}\wedge \underline{u'_i(2)}\wedge \underline{v_j}=\underline{0}$.
Hence~$\underline{u_i}\wedge \underline{u'_i(2)}\wedge\underline{\bar{v}_j}=0$, and~$\underline{\tilde{u}_i}\wedge\underline{\bar{v}_j}=0$.
Now we have proved
\begin{equation*}
\underline{\tilde{u}_i}\wedge\underline{\bar{v}_j}\begin{cases}
\neq 0, & \textrm{if~} i=j, \\
=0, & \textrm{if~} i<j,
\end{cases}
\end{equation*}
and it is enough to show that~$\underline{\tilde{u}_1},\ldots,\underline{\tilde{u}_m}$ are linear independent.
Suppose not, there exists not all zero~$c_1,\ldots,c_m$ such that~$c_1\underline{\tilde{u}_1}+\cdots+c_m\underline{\tilde{u}_m}=0$.
Suppose~$k$ is the maximum index such that~$c_k\neq 0$.
Then we have~$c_k\underline{\tilde{u}_k}=-(c_1\underline{\tilde{u}_1}+\cdots+c_{k-1}\underline{\tilde{u}_{k-1}})$,
and~$c_k\underline{\tilde{u}_k}\wedge\underline{\bar{v}_k}=-(c_1\underline{\tilde{u}_1}+\cdots+c_{k-1}\underline{\tilde{u}_{k-1}})\wedge\underline{\bar{v}_k}$.
Note that the left hand side is nonzero, but the right hand side is zero, so we get a contradiction.

For any~$i\in[m]$ and~$k\in[a_0]$,~$\underline{\alpha^{(i)}_k}$ can be linearly expressed by elements in~$\{\underline{\epsilon_1},\ldots,\underline{\epsilon_l}\}$,
so~$\underline{u_i}$ can be expressed in the form of
\begin{equation*}
\sum\limits_{1\leqslant s_1<\cdots<s_{a_0}\leqslant l} D_{s_1,\ldots,s_{a_0}}\underline{\epsilon_{s_1}}\wedge\cdots\wedge\underline{\epsilon_{s_{a_0}}},
\end{equation*}
where~$D_{s_1,\ldots,s_{a_0}}$ is a real number. Let
\begin{equation*}
P=\bigcup_{k=0}^{l-a_0-b_0}P_k,
\end{equation*}
where
\begin{equation*}
P_k=\{\underline{\epsilon_{s_1}}\wedge\cdots\wedge\underline{\epsilon_{s_k}}\wedge\underline{\gamma_{t_1}}\wedge\cdots\wedge\underline{\gamma_{t_{a-a_0-k}}}\mid
s_1<\cdots<s_k,t_1<\cdots<t_{a-a_0-k}\}.
\end{equation*}
Then~$\underline{u'_i(1)}$ can be linear expressed by elements in~$P$, and so~$\underline{\tilde{u}_i}$ can be expressed by elements in
\begin{equation*}
Q:=\bigcup_{k=0}^{l-a_0-b_0}Q_k,
\end{equation*}
where
\begin{equation*}
Q_k=\{\underline{\epsilon_{s_1}}\wedge\cdots\wedge\underline{\epsilon_{s_{a_0+k}}}\wedge\underline{\gamma_{t_1}}\wedge\cdots\wedge\underline{\gamma_{t_{a-a_0-k}}}\mid
s_1<\cdots<s_{a_0+k},t_1<\cdots<t_{a-a_0-k}\}.
\end{equation*}
Then we have
\begin{equation*}
m\leqslant|Q|=\sum_{k=0}^{l-a_0-b_0}|Q_k|=\sum_{k=0}^{l-a_0-b_0}\binom{l}{a_0+k}\binom{a+b-l}{a-a_0-k}.
\end{equation*}

\textbf{Case 2:}~$n+l>a+b$.

Using Lemma~\ref{GP3}, there exists a linear mapping~$\phi:V\rightarrow \mathbb{R}^{a+b}$ that keep the dimension of
\begin{equation*}
A_i,\quad B_i,\quad W,\quad W+A_i,\quad W+B_i,\quad A_i+B_j,\quad i,j\in[m].
\end{equation*}
Here the condition~$l\leqslant\min\{a_0+b,a+b_0\}$ guarantees that the dimension of all of these spaces are no more than~$a+b$.
Note that
\begin{equation*}
\begin{array}{c}
  \dim(\phi(W)+\phi(A_i))=\dim(\phi(W))+\dim(\phi(A_i))-\dim(\phi(W)\cap\phi(A_i)), \\
  \dim(W+A_i)=\dim(W)+\dim(A_i)-\dim(W\cap A_i), \\
  \dim(\phi(W))=\dim(W),\quad\dim(\phi(A_i))=\dim(A_i), \\
  \dim(\phi(W)+\phi(A_i))=\dim(\phi(W+A_i))=\dim(W+A_i),
\end{array}
\end{equation*}
we have
\begin{equation*}
\dim(\phi(W)\cap\phi(A_i))=\dim(W\cap A_i).
\end{equation*}
Similarly, we have
\begin{equation*}
\dim(\phi(W)\cap\phi(B_i))=\dim(W\cap B_i),\quad\dim(\phi(A_i)\cap\phi(B_j))=\dim(A_i\cap B_j).
\end{equation*}
In this way, we can replace~$V,W,A_i,B_i$ to~$\mathbb{R}^{a+b},\phi(W),\phi(A_i),\phi(B_i)$, and they still satisfy the condition of the theorem.
So we reduced the general case to the case~$n+l=a+b$ that we have already proved.

\section{Two examples}

In this section, two examples of skew~\bol~system will be constructed,
to show the tightness of Inequality~(\ref{Eq:main}), and the necessity of the condition that~$l$ is small.

Suppose~$V=\mathbb{R}^{n+l}$, and~$W$ is an~$l$-dimensional subspace of~$V$ for some nonnegative~$n,l$.
Note that in the two examples,~$n,l$ are valued differently.
Pick~$\underline{\epsilon_1},\ldots,\underline{\epsilon_l}$, a basis of~$W$,
and expand it to~$\underline{\epsilon_1},\ldots,\underline{\epsilon_l},\underline{\gamma_1},\ldots,\underline{\gamma_n}$, a basis of~$V$.

\textbf{Example 1.}
Suppose~$a\geqslant a_0\geqslant 0$,~$b\geqslant b_0\geqslant 0$, and~$l\leqslant\min\{a_0+b,a+b_0\}$, let~$n=a+b-l$.
For a fixed~$k\in\{0,1,\ldots,l-a_0-b_0\}$, and~$S\in\binom{[l]}{a_0+k},T\in\binom{[n]}{a-a_0-k}$, denote
\begin{gather*}
    S=\{s_1,\ldots,s_{a_0+k}\}\subseteq[l],\quad T=\{t_1,\ldots,t_{a-a_0-k}\}\subseteq[n], \\
    P=\{p_1,\ldots,p_{l-a_0-k}\}:=[l]\setminus S,\quad Q=\{q_1,\ldots,q_{n-a+a_0+k}\}:=[n]\setminus T.
\end{gather*}
Let
\begin{gather*}
    A(S,T)=\spa\{\underline{\epsilon_{s_1}},\ldots,\underline{\epsilon_{s_{a_0}}},
\underline{\gamma_{q_1}}+\underline{\epsilon_{a_0+1}},\ldots,\underline{\gamma_{q_k}}+\underline{\epsilon_{s_{a_0+k}}},
\underline{\gamma_{t_1}},\ldots,\underline{\gamma_{t_{a-a_0-k}}}\}, \\
    B(S,T)=\spa\{\underline{\epsilon_{p_1}},\ldots,\underline{\epsilon_{p_{b_0}}},
\underline{\gamma_{t_1}}+\underline{\epsilon_{p_{b_0+1}}},\ldots,\underline{\gamma_{t_{l-a_0-b_0-k}}}+\underline{\epsilon_{p_{l-a_0-k}}},
\underline{\gamma_{q_1}},\ldots,\underline{\gamma_{q_{b+a_0-l+k}}}\}.
\end{gather*}
Then by the definition,~$A(S,T)$ and~$B(S,T)$ are~$(a,a_0)$-type and~$(b,b_0)$-type subspaces correspondingly, and~$\dim(A(S,T)\cap B(S,T))=0$ for any~$S$ and~$T$.

Let~$\mathcal{P}_k$ be the families of all pairs~$(A(S,T),B(S,T))$ defined above (for the fixed~$k$), and
\begin{equation*}
\mathcal{P}=\bigcup_{k=0}^{l-a_0-b_0}\mathcal{P}_k.
\end{equation*}
Note that for fixed~$k$ there are~$\binom{l}{a_0+k}$ choices for~$S$, and~$\binom{a+b-l}{a-a_0-k}$ choices for~$T$, so
\begin{equation*}
|\mathcal{P}|=\sum_{k=0}^{l-a_0-b_0}\binom{l}{a_0+k}\binom{a+b-l}{a-a_0-k}.
\end{equation*}

Index these~$A(S,T)$'s and~$B(S,T)$'s such that
\begin{itemize}
  \item $A(S,T)$ and~$B(S,T)$ have the same index for any~$S,T$;
  \item $(A(S,T),B(S,T))\in\mathcal{P}_k$ has a smaller index than~$(A(S',T'),B(S',T'))\in\mathcal{P}_{k'}$ if~$k<k'$.
\end{itemize}

To prove~$\mathcal{P}$ is a skew~\bol~system, we need to check~$\dim(A_i\cap B_j)>0$ if~$i<j$. Suppose~$A_i=A(S,T),(A_i,B_i)\in\mathcal{P}_k$ and~$B_j=B(S',T'),(A_j,B_j)\in\mathcal{P}_{k'}$.
Let~$P=[l]\setminus S, Q=[n]\setminus T, P'=[l]\setminus S', Q'=[n]\setminus T'$.
By the principle of index,~$i<j$ implies~$k\leqslant k'$. There are~$3$ possibilities (as shown in the diagram) to discuss.

\textbf{Case 1:}~$k<k'$. Since~$|T|+|Q'|>n$, there exists a~$c\in T\cap Q'$.

\textbf{Case 2:}~$k=k'$ and~$T\neq T'$. Again there is a~$c\in T\cap Q'$. For both of these two cases, we have~$\underline{\gamma_c}\in A(S,T)\cap B(S',T')$,
so~$A_i\cap B_j\neq\emptyset$.

\textbf{Case 3:}~$k=k'$ and~$T=T'$. Then~$S\neq S'$, so~$S\cap P'\neq\emptyset$. Suppose~$c\in S\cap P'$.
Then~$\underline{\epsilon_c}+\underline{\mu}\in A(S,T)$,
where either~$\underline{\mu}=\underline{0}$ or~$\underline{\mu}=\underline{\gamma_q}$ for some~$q\in Q$.
Similarly,~$\underline{\epsilon_c}+\underline{\nu}\in B(S',T')$,
where either~$\underline{\nu}=\underline{0}$ or~$\underline{\nu}=\underline{\gamma_t}$ for some~$t\in T'$.
Using the fact~$T=T',Q=Q'$, we have~$\underline{0}\neq\underline{\epsilon_c}+\underline{\mu}+\underline{\nu}\in A(S,T)\cap B(S',T')$,
so~$A_i\cap B_j\neq\emptyset$.

\begin{table}[h]
\centering
\setlength{\tabcolsep}{2mm}{
\begin{tabular}{|c|c|c|c|c|c|c|c|c|c|c|c|cl}
  \multicolumn{6}{|c|}{$\epsilon_1,\ldots,\epsilon_l$}   &  \multicolumn{6}{|c|}{$\gamma_1,\ldots,\gamma_n$}  & &\\
  \cline{1-12}
  \multicolumn{1}{c}{\phantom{a}} & \multicolumn{1}{c}{\phantom{a}} & \multicolumn{1}{c}{\phantom{a}} & \multicolumn{1}{c}{\phantom{a}} & \multicolumn{1}{c}{\phantom{a}}
   & \multicolumn{1}{c}{\phantom{a}} & \multicolumn{1}{c}{\phantom{a}} & \multicolumn{1}{c}{\phantom{a}}
    & \multicolumn{1}{c}{\phantom{a}} & \multicolumn{1}{c}{\phantom{a}} & \multicolumn{1}{c}{\phantom{a}} & \multicolumn{1}{c}{\phantom{a}} & &\\
  \cline{1-12}
  \multicolumn{3}{|c|}{$S$} & \multicolumn{3}{|c|}{$P$} & \multicolumn{3}{|c|}{$T$} & \multicolumn{3}{|c|}{$Q$} & & \\
  \cline{1-12}
  \multicolumn{14}{c}{} \\
  \cline{7-12}
  \multicolumn{6}{c}{} & \multicolumn{2}{|c|}{$T'$}  & \multicolumn{4}{|c|}{$Q'$} & & Case 1:~$k<k'. $ \\
  \cline{7-12}
  \multicolumn{14}{c}{} \\
  \cline{7-12}
  \multicolumn{6}{c}{} & \multicolumn{2}{|c|}{$T'$} & $Q'$ & $T'$ & \multicolumn{2}{|c|}{$Q'$} & & Case 2:~$k=k', T\neq T'. $\\
  \cline{7-12}
  \multicolumn{14}{c}{} \\
  \cline{1-12}
  \multicolumn{2}{|c|}{$S'$} & $P'$ & $S'$ & \multicolumn{2}{|c|}{$P'$} & \multicolumn{3}{|c|}{$T'=T$} & \multicolumn{3}{|c|}{$Q'=Q$} & & Case 3:~$k=k', T=T', S\neq S'$. \\
  \cline{1-12}
\end{tabular}}
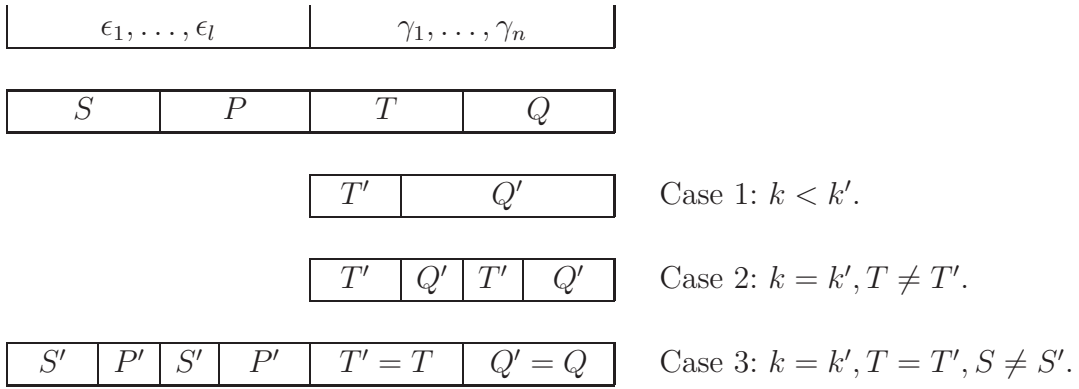
\caption*{A schematic diagram for the three cases}
\end{table}

Then the exact value of~$M(a,a_0,b,b_0,n,l)$ is determined for the case~$l$ is not too big by Theorem~\ref{Th:Spaces} and the example above.
Now we are going to discuss the case~$l>\min\{a_0+b,a+b_0\}$. Without lose of generality, we may assume that~$a+b_0\leqslant a_0+b$.
Using a similar tactic, we can construct a skew~\bol~system of subspaces~$\mathcal{P}$ with as large cardinality as possible.

\textbf{Example 2.}
For any fixed~$k\in\{a_0,\ldots,a\}$, and~$S\in\binom{[a+b_0]}{k},T\in\binom{[b-b_0]}{a-k}$, denote
\begin{gather*}
    S=\{s_1,\ldots,s_k\}\subseteq[a+b_0],\quad T=\{t_1,\ldots,t_{a-k}\}\subseteq[b-b_0], \\
    P=\{p_1,\ldots,p_{a+b_0-k}\}:=[a+b_0]\setminus S,\quad Q=\{q_1,\ldots,q_{b-b_0-a+k}\}:=[b-b_0]\setminus T.
\end{gather*}
Let
\begin{gather*}
    A(S,T)=\spa\{\underline{\epsilon_{s_1}},\ldots,\underline{\epsilon_{s_{a_0}}},
\underline{\gamma_{q_1}}+\underline{\epsilon_{a_0+1}},\ldots,\underline{\gamma_{q_{k-a_0}}}+\underline{\epsilon_{s_k}},
\underline{\gamma_{t_1}},\ldots,\underline{\gamma_{t_{a-k}}}\}, \\
    B(S,T)=\spa\{\underline{\epsilon_{p_1}},\ldots,\underline{\epsilon_{p_{b_0}}},
\underline{\gamma_{t_1}}+\underline{\epsilon_{p_{b_0+1}}},\ldots,\underline{\gamma_{t_{a-k}}}+\underline{\epsilon_{p_{a+b_0-k}}},
\underline{\gamma_{q_1}},\ldots,\underline{\gamma_{q_{b-b_0-a+k}}}\}.
\end{gather*}
Then by the definition,~$A(S,T)$ and~$B(S,T)$ are~$(a,a_0)$-type and~$(b,b_0)$-type subspaces correspondingly, and~$\dim(A(S,T)\cap B(S,T))=0$ for any~$S$ and~$T$.

Let~$\mathcal{P}_k$ be the families of pairs~$(A(S,T),B(S,T))$, and
\begin{equation*}
\mathcal{P}=\bigcup_{k=a_0}^a\mathcal{P}_k.
\end{equation*}
Then we have
\begin{equation*}
|\mathcal{P}|=\sum_{k=a_0}^a\binom{a+b_0}{k}\binom{b-b_0}{a-k}.
\end{equation*}

Index these~$A(S,T)$'s and~$B(S,T)$'s such that
\begin{itemize}
  \item $A(S,T)$ and~$B(S,T)$ have the same index for any~$S,T$;
  \item $(A(S,T),B(S,T))\in\mathcal{P}_k$ has a smaller index than~$(A(S',T'),B(S',T'))\in\mathcal{P}_{k'}$ if~$k<k'$.
\end{itemize}
Similar with the discussion in the former example, the families~$\mathcal{P}$ is a skew~\bol~system.

Now we find a lower bound for the maximum cardinality~$M(a,a_0,b,b_0,n,l)$ for the case~$l>\min\{a_0+b,a+b_0\}$ and~$a+b_0\leqslant a_0+b$.
Note that~$\binom{a+b}{a}$ is always a natural upper bound, so in this case we have
\begin{equation}\label{Eq:lBig}
\sum_{k=a_0}^a\binom{a+b_0}{k}\binom{b-b_0}{a-k}\leqslant M(a,a_0,b,b_0,n,l)\leqslant\binom{a+b}{a}=\sum_{k=0}^a\binom{a+b_0}{k}\binom{b-b_0}{a-k}.
\end{equation}
The difference between the two sides of~(\ref{Eq:lBig}) are the terms of~$k<a_0$.
Specifically, if~$a_0+b=a+b_0$ then~$\binom{b-b_0}{a-k}=0$ for any~$k<a_0$, so these terms vanished. Then we have the following statement.

\begin{mingti}\label{Th:lBig}
If~$l>a_0+b=a+b_0$ then~$M(a,a_0,b,b_0,n,l)=\binom{a+b}{a}$.
\end{mingti}

Note that the right hand side of~(\ref{Eq:main}) is strictly smaller than~$\binom{a+b}{a}$,
so the discussion above shows that the condition~$l\leqslant\min\{a_0+b,a+b_0\}$ in Theorem~\ref{Th:Spaces} is not redundant.

Similar with~(\ref{Eq:lBig}), we can also prove the following lower bound for the case~$l>\min\{a_0+b,a+b_0\}$ and~$a+b_0> a_0+b$.
\begin{equation*}
M(a,a_0,b,b_0,n,l)\geqslant\sum_{k=a_0}^a\binom{a_0+b}{k}\binom{a-a_0}{a-k}.
\end{equation*}

\begin{center}
\textbf{Acknowledgments}
\end{center}

K. Wang is supported by the National Key R\&D Program of China (No.2020YFA0712900) and NSFC (12071039,12131011).

\bibliography{bib}

\end{document}